\documentclass[10pt,leqno,twoside]{article}

\usepackage[latin2]{inputenc}
\usepackage{amsmath}
\usepackage{amsfonts}
\usepackage{amssymb}
\usepackage{makeidx}

\newtheorem{statement}{Statement}
\newtheorem{theorem}{Theorem}

\newtheorem{corollary}{Corollary}
\newtheorem{lemma}{Lemma}

\newfont{\eurorm}{eurm10 scaled 1100}
\newfont{\eurosm}{eurm10 scaled 800}

\newfont{\newit}{cmfi10 scaled 1100}
\newfont{\newsit}{cmfi10 scaled 1000}
\newfont{\newsmit}{cmfi10 scaled 800}
\newfont{\newssmit}{cmfi10 scaled 600}

\newfont{\cirilrm}{wncyr10 scaled 1200}
\newfont{\cirilbf}{wncyb10 scaled 1200}
\newfont{\cirilsf}{wncyss10 scaled 1200}
\newfont{\cirilit}{wncyi10 scaled 1200}
\newfont{\cirilsc}{wncysc10 scaled 1200}

\def\mm#1{\mbox{\eurorm #1}}
\def\mmsm#1{\mbox{\eurosm #1}}

\begin{document}

$\,$

\bigskip

\centerline{\Large \bf Refined estimates  and generalizations of
inequalities}

\smallskip

\centerline{\Large \bf related to the arctangent function and
Shafer's inequality}

\smallskip

\begin{center}
{\em Branko Male\v sevi\' c, Marija Ra\v sajski, Tatjana Lutovac}
\end{center}

\begin{center}
University of Belgrade, School  of Electrical Engineering,   \\[0.0 ex]
Department of Applied Mathematics, Serbia
\end{center}

\medskip
\noindent
{\small \textbf{Abstract.} In this paper we give some sharper refinements and generalizations of inequalities
related to {\sc Shafer}'s inequality for the arctangent function, stated in Theorems 1, 2 and 4  in \cite{MS_2014},
by {\sc C. Mortici } and {\sc H.$\,$M. Srivastava}.}


\footnote{$\!\!\!\!\!\!\!\!\!\!\!\!\!\!$ \scriptsize Emails: Branko
Male\v sevi\' c {\tt $<$branko.malesevic@etf.rs$>$}, Marija Ra\v
sajski {\tt $<$marija.rasajski@etf.rs$>$}, Tatjana Lutovac {\tt
$<$tatjana.lutovac@etf.rs$>$}}

\vspace*{-2.5 mm}

{\small \tt MSC: Primary 33B10; Secondary 26D05}

\section{Introduction}

\medskip

Inverse trigonometric functions play an important role and have many applications in engineering
\cite{deAbreu_2009}, \cite{deAbreu_2012}, \cite{Alirezaei_Mathar_2014} and \cite{Bercu_2017}.
In particular, the arctangent function and various related inequalities have been studied and
effectively applied to problems in fundamental sciences and many areas of engineering,
such as electronics, mechanics, and aeronautics \cite{Alirezaei_Mathar_2014}, \cite{Alirezaei_2013}, \cite{Alirezaei_2014};
see also \cite{CloudDrachmanLebedev_2014}.

\smallskip
Various approximations of the arctangent function can be found in \cite{MS_2014},
\cite{Shafer_1966}-\cite{Nishizawa_2017}; see also \cite{DSM_1970}, \cite{GVM_2014}.
One of the inequalities that attracted attention of many authors is {\sc Shafer}'s inequality \cite{Shafer_1966}:
\begin{equation}\label{Shafer_ineq_66}\displaystyle
\frac{3x}{1+2 \sqrt{\,1+x^2\, }\,}  \, <  \,  \mbox{\rm arctan} \,
x,
\end{equation}
which holds for $x > 0$; see also \cite{Shafer_1974}-\cite{Shafer_1978}.

\smallskip
Recently, in \cite{MS_2014}, {\sc Mortici } and {\sc Srivastava}
proved  the following results, cited here as Statements 1, 2 and 3,
related to the above inequality. These results are the starting
point of our research.

\begin{statement}
\label{MS_statement_1} $($Theorem $1$, {\rm
\cite{MS_2014}}$)$ For every  $x > 0$, the following
two-sided inequality holds$\,:$
\begin{equation}\label{MS_1} \displaystyle
\frac{3x}{1 + 2\sqrt{1 + x^2}} +   a(x)  < \mbox{\rm arctan}\,x  <
\frac{3x}{1 + 2\sqrt{1 + x^2}} + b(x),
\end{equation}
where $\displaystyle a(x) = \mbox{\small  $\displaystyle\frac{1}{180}x^{5} - \frac{13}{1512}x^{7}$}$
and $b(x)= \mbox{\small $\displaystyle \frac{1}{180}x^{5}$}$.
\end{statement}

\medskip

\begin{statement}
\label{MS_statement_2} $($Theorem $2$, {\rm
\cite{MS_2014}}$)$ For every  $x > 0$, it is
asserted that
\begin{equation}\label{MS_2} \displaystyle
\frac{3x+ c(x)}{1 + 2\sqrt{1 + x^2}}   < \mbox{\rm arctan}\,x <
\frac{3x+ d(x)}{1 + 2\sqrt{1 + x^2}},
\end{equation}
where $c(x) = \mbox{\small $\displaystyle \frac{1}{60}x^{5} -
\frac{17}{840}x^{7}$}$ and $d(x)= \mbox{\small $\displaystyle
\frac{1}{60}x^{5}$}$.
\end{statement}

\medskip

\begin{statement}\label{MS_statement_4} $($Theorem $4$, {\rm
\cite{MS_2014}}$)$ For every  $x > 0$, it is
asserted that
\begin{equation}\label{MS_4}  \displaystyle
-\frac{1}{12}x^{3} \, < \mbox{\rm arctan}\,x   -    \frac{2x}{1 +
\sqrt{1 + x^2}} <  \, -\frac{1}{12}x^{3} \, + \, \frac{3}{40}x^{5}.
\end{equation}
\end{statement}

The main results of this paper are refined estimates and
generalizations  of the inequalities  given  in Statements 1, 2 and
3.  Although  the inequalities (\ref{MS_1}), (\ref{MS_2}) and
(\ref{MS_4})  hold for $x>0$,
considering them
in a neighborhood of zero is of primary importance, as noted in
\cite{MS_2014}.

\section{Main results}


First, let us recall some well-known power series expansions that
will be used in our proofs.

\smallskip
For $|x| \leq 1$,
\begin{equation}
\label{label_A} \mbox{\rm arctan}\,x \,=
\sum\limits_{m=0}^{\infty}{(-1)^{m} \,A(m) \, x^{2m+1}},
\end{equation}
where
\begin{equation}
A(m)=\displaystyle\frac{1}{2m+1}.
\end{equation}

 For $|x| \leq 1$,
\begin{equation}
\label{label_K} \sqrt{1+x^2} \,= 1 + \sum_{m=0}^{\infty}{(-1)^{m}
\,K(m) \, x^{2m+2}},
\end{equation}
where
\begin{equation}
K(m) = \displaystyle\frac{(2m)!}{m!\, (m+1)!\,2^{2m+1}}.
\end{equation}

\smallskip
 The following power series expansion holds:
\begin{equation}
\displaystyle\frac{3x}{1+2\sqrt{1+x^2}}
=
\sum_{m=0}^{\infty}{B(m) \, x^{2m+1}},
\end{equation}
where $|x|\!\leq\!\mbox{\small $\displaystyle\frac{\sqrt{3}}{2}$}$,
with $B(0)=1$, $B(1)=-\mbox{\small $\displaystyle\frac{1}{3}$}$, and
for $m\geq2$:
\begin{equation}
\label{S(m)_expl}
B(m)
=
\frac{(-1)^{m}4^{m-1}}{3^{m}}
\left(
1-8\sum\limits_{i=2}^{m}{\frac{(2i-2)!}{(i-1)!\,i!\,2^{2i-1}}\!\left(\frac{3}{4}\right)^{\!\!i\,}}
\right).
\end{equation}
Power series coefficients are calculated by applying {\sc c}'s
product to the power series expansions arising from the following
transformation of the corresponding function:
\begin{equation}
\frac{3x}{1+2\sqrt{1+x^2}}
= - \,
\frac{x\left(1 -  2 \sqrt{1+x^2}\,\right)}{1 + \mbox{\small $\displaystyle\frac{4}{3}$}x^2}.
\end{equation}

It is easy to prove that  sequence $\left\{ {B\left( m \right)}
\right\}_{m \in N_{0}}$  for $m \geq 1$ satisfies the recurrence
equation:
\begin{equation}
\label{Difference_Eq_B} \begin{array}{lcl}\displaystyle
B(m+1)+\frac{4}{3}B(m)
\!&\!\!=\!\!&\!
\displaystyle (-1)^m \frac{(2\,m \,- \, 1)!}{\, 2^{2m-1} \, (m \, - \, 1)! \, (m\, + \, 1)! \,}  \\[1.5em]
\!&\!\!=\!\!&\!
\displaystyle {2 \, ( - 1)^m}\, \frac{{(2m - 1)!!}}{{(2m + 2)!!}}.
\end{array}
\end{equation}

\bigskip
\subsection{   Refinements of the inequalities in Statement 1}  

Before we proceed to Theorem 1, which represents an improvement and
generalization of Statement \ref{MS_statement_1}, we need the
following lemmas.
\medskip
\begin{lemma}$\:$\\
Let $\beta(0)\!=\!1$,  $\beta(1)\!=\!-\mbox{\small $\displaystyle\frac{1}{3}$}$
and $\beta(m)\!=\!\displaystyle\frac{3(-1)^m}{2^{2m+1}}\!\sum_{k=0}^{\infty}{\frac{(2k+2m-1)!}{(k+m-1)!\,(k+m+1)!\,}
\!\left(\mbox{\small $\displaystyle\frac{3}{16}$}\right)^{\!k}}\!,$ for $m\!\geq\!2$.
The sequence $\left\{\beta(m)\right\}_{m \in N_{0}}$ for $m\!\geq\!1$
satisfies the recurrence relation \hfill

\smallskip \noindent {\rm (\ref{Difference_Eq_B})}.
\end{lemma}

\smallskip
\noindent
{\bf Proof.} In the proof of this lemma we use the  {\sc Wilf}-{\sc Zeilberger}
method \cite{PWZ_1996}, \cite{T_2010} and \cite{PZ_2013}.
(The~same approach we used in \cite{MRL_2017}.)

\smallskip
The assertion   is   obviously  true for   $m=1$.

\smallskip
Let $m \geq 2$ and
$$g(k, \,m)
=
{\displaystyle
\frac {(2\,k + 2\,m - 1)\mathrm{!}}{(k + m - 1)\mathrm{!}\,(k + m +1)\mathrm{!}}}
\quad \mbox{and} \quad
\phi(m) = {\displaystyle \frac {3\,(-1)^{m}}{2^{2\,m+1}}}.
$$
Then we have
\[
\mathit{\beta(m)} =  \phi (m) {\displaystyle \sum
_{k=0}^{\infty }} \,\mathrm{g}(k, \,m) \mbox{\small $\left(\!{\displaystyle \frac {3}{16}}\!\right)$}^{k} \! .
\]
Further we have

\smallskip
\begin{equation}\label{relacija-4}
\begin{array}{rcl}
\mathit{\beta}(m + 1) + {\displaystyle \frac {4}{3}} \,\mathit{\beta}(m)
\!&\!\!=\!\!&\! {\displaystyle \sum_{k=0}^{\infty }} \,{\displaystyle \left( \phi(m+1) g(k, \,m+1)
+
{\displaystyle \frac{4}{3}} \phi(m) g(k, \,m) \right)\!\!\left(\mbox{\small $\displaystyle\frac{3}{16}$}\right)^{\!k}}  \\
\!&\!\!=\!\!&\!
{\displaystyle \frac{(-1)^{m}}{2^{(2\,m - 1)}}} \cdot {\displaystyle \sum _{k=0}^{\infty }}
{\displaystyle\frac{(2\,k + 2\,m + 13)\,(2 \,k + 2\,m - 1)\mathrm{!}\,}{ (8\,k + 8\,m +
16)\,(k + m - 1)\mathrm{!}\,(k + m + 1)\mathrm{!}}} \!\left(\mbox{\small $\displaystyle\frac {3}{16}$}\right)^{\!k}.
\end{array}
\end{equation}
Consider now the sequence $\left\{ \mm{S}(m)\right\}_{m \in N, \,
m\geq 2}$,  where

\smallskip
\begin{equation}\label{S(m)} \displaystyle
\mm{S}(m) \!=\!\! \displaystyle\sum_{k=0}^{m-1}\mm{F}(m,k)
\quad \mbox{and} \quad
\displaystyle \mm{F}(m,k) = \mbox{\small $
{\displaystyle \frac{(2m\!+\!2k\!+\!13) \,(2m\!+\!2k\!-\!1)\mathrm{!}}{(8m\!+\!8k\!+\!16)\,(m\!+\!k\!-\!1)\mathrm{!}\,(m\!+\!k\!+\!1)\mathrm{!}}}$}
\!\left({\mbox{\small $\displaystyle\frac{3}{16}$}}\right)^{\!k}\!.
\end{equation}
Consider the function${}^{\ast)}$\footnote{$\!\!\!\!\!\!\!\!\!\!\!{}^{\ast)}\,$An
algorithm  for determining function $\mmsm{G}(m,k)$ for a given function $\mmsm{F}(m,k)$ is described in \cite{PWZ_1996}.$\;$Note 
that the pair of discrete functions ${\big (} \mmsm{F}(m,k), \, \mmsm{G}(m,k){\big )}$ is the so-called {\sc Wilf}-{\sc Zeilberger} pair.}
$$
\mm{G}(m,k) = \mbox{\small $\displaystyle\frac{( - 8\,k - 8\,m - 16)\,(2\,k + 2\,m - 1)
\mathrm{!}}{(8\,k + 8\,m + 16)\,(k + m - 1)\mathrm{!}\,(k + m + 1)\mathrm{!}}$}
\!\left(\mbox{\small $\displaystyle\frac{3}{16}$}\right)^{\!k}\!,
$$
where  $m \!\in\! N_{0}$ and $k \!\in\! N_{0}$. It is not hard to
verify that functions $\mm{F}(m,k)$ and $\mm{G}(m,k)$ satisfy the
following relation:
\begin{equation}
\label{F-G}
\mm{F}(m, \,k)=\mm{G}(m, \,k + 1) - \mm{G}(m, \,k).
\end{equation}
If we sum both sides of (\ref{F-G}) over all $k\in N_{0}$,
 we get the following relation:
$$\begin{array}{l}
\mm{S}(m) = - \mm{G}(m,0).
\end{array}
$$
Finally,  as
$$\label{initial-conditions}
\mm{G}(m,0)= {\displaystyle - \frac {(2\,m - 1)\mathrm{!}}{(m -
1)\mathrm{!} \,(m + 1 )\mathrm{!}}}
$$
we have \begin{equation}\label{suma} \mm{S}(m)= {\displaystyle \frac
{(2\,m - 1)\mathrm{!}}{(m - 1)\mathrm{!}\,(m + 1 )\mathrm{!}}}.
\end{equation}
Therefore from (\ref{relacija-4}) and (\ref{suma})  we conclude that
$$
\begin{array}{rcl}
\mathit{\beta}(m + 1) + {\displaystyle \frac {4}{3}} \,\mathit{\beta}(m)
\!&\!\!=\!\!&\!
{\displaystyle \frac {(-1)^{m}}{2^{2\,m - 1}}} \cdot \mm{S}(m) \\[1.25em]
\!&\!\!=\!\!&\!
{\displaystyle \frac {(-1)^{m}}{2^{2\,m - 1}}} \cdot {\displaystyle \frac {(2\,m - 1)\mathrm{!}}{(m - 1)\mathrm{!}\,(m + 1 )\mathrm{!}}}.
\end{array}
$$

\hfill $\Box$

\break

\begin{corollary} 
Given that  the sequences $\left\{\mm{B}(m) \right\}_{m \in N_{0}}$ and $\left\{\mathit{\beta}(m) \right\}_{m \in N_{0}}$ 
satisfy the same recurrence relation and as they agree for $m=0$ and  $m=1$, we conclude that
\begin{equation}\label{jednakost}
\mm{B}(m) =\mathit{\beta}(m),~~~ \mbox{for} ~\, m\!\in\!N_{0}.
\end{equation}
\end{corollary} 




\medskip
Let us introduce the notation:
\begin{equation}
C(m)=A(m)-B(m),
\end{equation}
where $C(0)=C(1)=0$, and for $m\geq2$ the following holds:
\begin{equation}
C(m) \!=\! \displaystyle\frac{1}{2m+1} -
\displaystyle\frac{4^{m-1}}{3^{m}} \left(
1-8\sum\limits_{i=2}^{m}{\frac{(2i-2)!}{(i-1)!\,i!\,2^{2i-1}}\!\left(\frac{3}{4}\right)^{\!\!i\,}}
\right).
\end{equation}

Thus, we have the power series expansion:
\begin{equation}
f(x) = \mbox{arctan}\,x - \displaystyle\frac{3x}{2+\sqrt{1+x^2}} =
\sum_{m=0}^{\infty}{{(-1)^{m}}C(m) x^{2m+1}},
\end{equation}
for $m \!\in\! N_0$ and $|x| \!\leq\! \mbox{\small
$\displaystyle\frac{\sqrt{3}}{2}$}$.

\medskip
Let us introduce the notation:
\[
\beta^{+}(m)=
\displaystyle\frac{3}{2^{2m+1}}
\sum_{k=0}^{\infty}{\frac{(2k+2m-1)!}{
(k+m-1)! \, (k+m+1)! \, 2^{2k}} \! \left(\mbox{\small $\displaystyle\frac{3}{4}$}\right)^{\!k}}, m \!\in\! N_0
\]
\[
{\beta_1}(m)
=
\mathop \sum \limits_{k = 0}^{+\infty}
\frac{{\mathop \prod \nolimits_{n = 0}^{k - 1} \left( {m + \frac{1}{2} + n} \right)}}{
{\mathop \prod \nolimits_{n = 0}^{k - 1} \left( {m + 2 + n} \right)}}{\left( {\frac{3}{4}} \right)^k}\!,
\;\;\beta_1^k\left( m \right) = \frac{{\mathop \prod \nolimits_{n = 0}^{k - 1} (m + \frac{1}{2} + n)}}{
{\mathop \prod \nolimits_{n = 0}^{k - 1} (m + 2 + n)}}, \, m \!\in\! N_0.
\]
\begin{lemma}$\:$\\[0.5 em] 
For $m \in N_0$ the following  holds:
\[\beta^{+}(m) = \frac{3}{2}\frac{{(2m - 1)!!}}{{(2m + 2)!!}}{\beta_1}(m).\]
\end{lemma}
{\bf Proof.}
{\small 
\[
\!\!
\begin{array}{lcl} 
\displaystyle\frac{3}{{{2^{2m + 1}}}}\frac{{(2m \!+\! 2k \!-\! 1)!}}{{(m \!+\! k \!-\! 1)!(m \!+\! k \!+\! 1)!{2^{2k}}}}
\!&\!\!=\!\!&\! \displaystyle \frac{3}{2}\frac{1}{{(2m \!+\! 2)!!}}\frac{{(2m \!+\! 2k \!-\! 1)!!}}{{(m \!+\! k \!+\! 1)!}}\frac{(m \!+\! 1)!}{{{2^k}}}  \\ [1.25em]
\!&\!\!=\!\!&\! \displaystyle \frac{3}{2}\frac{{(2m \!-\! 1)!!}}{{(2m \!+\! 2)!!}}\frac{{(2m \!+\! 1)(2m \!+\! 3) \ldots (2m \!+\! 2k \!-\! 1)}}{{(m \!+\! 1)!(m \!+\! 2)(m \!+\! 3) \ldots (m \!+\! k \!+\! 1)}}\frac{(m \!+\! 1)!}{{{2^k}}}  \\ [1.25em]
\!&\!\!=\!\!&\! \displaystyle \frac{3}{2}\frac{{(2m \!-\! 1)!!}}{{(2m \!+\! 2)!!}}\frac{{\mathop \prod \nolimits_{n = 0}^{k - 1} (m \!+\! \frac{1}{2} \!+\! n)}}{{\mathop \prod \nolimits_{n = 0}^{k-1}(m \!+\! 2 \!+\! n)}}.
\end{array}\]}
\hfill  $\Box$
\begin{lemma}$\:$\\[0.5 em] 
For $m \in N_0$ the following holds:
\[\frac{{8(m + 2)}}{{2m + 13}} = \mathop \sum \limits_{k = 0}^\infty  {\left( {\frac{{m + \frac{1}{2}}}{{m + 2}}} \right)^k}{\left( {\frac{3}{4}} \right)^k} < {\beta_1}(m) < \mathop \sum \limits_{k = 0}^\infty  {\left( {\frac{3}{4}} \right)^k} = 4.\]
\end{lemma}
{\bf Proof.} The statement immediately follows from the inequalities:
\[ {\left( {\frac{{m + \frac{1}{2}}}{{m + 2}}} \right)^k} \!< {\beta_1^k}(m) <  1.\]
\hfill $\Box$
\begin{lemma}$\:$\\[0.5 em] 
For $m \in N_0$ the following holds:
\[{\beta_1}(m + 1) > {\beta_1}(m)\frac{{m + 2}}{{m+\frac{1}{2}}}\left( {1 - \frac{3}{2}\frac{1}{{m + 2}}} \right).\]
\end{lemma}
{\bf Proof.}
\[\begin{array}{lcl}
\displaystyle \beta_1^k\left( {m + 1} \right)
\!&\!\!=\!\!&\!
\displaystyle \beta_1^k\left( m \right)\frac{{(m + 2)(m + \frac{1}{2} + k)}}{{(\frac{1}{2} + m)(m + 2 + k)}} \\[1.25em]
\!&\!\!=\!\!&\!
\displaystyle \beta_1^k\left( m \right)\frac{{m + 2}}{{m+\frac{1}{2} }}\left( {1 - \frac{3}{2}\frac{1}{{m + 2 + k}}} \right) \\[1.25em]
\!&\!\!>\!\!&\!
\displaystyle \beta_1^k\left( m \right)\frac{{m + 2}}{{ m+\frac{1}{2}}}\left( {1 - \displaystyle\frac{3}{2}\frac{1}{{m + 2}}} \right).
\end{array}\] 
\hfill  $\Box$
\begin{theorem}
\label{improvement-MS-th-1} 
For the real analytic function:
\begin{equation}
f(x)
=
\mbox{\rm arctan}\,x
-
\displaystyle\frac{3x}{2+\sqrt{1+x^2}}
\end{equation}\label{improvement1}
the following inequalities hold for $k \!\in\! N$ and $x \!\in\! {\Big (}0, \mbox{\small $\dfrac{\sqrt{3}}{2}$}{\Big ]}\!:$
\begin{equation}
\label{Mara22}
\sum\limits_{m = 0}^{2k + 1} {{{\left( { - 1} \right)}^m}C(m) {x^{2m
+ 1}}}  <f(x) < \sum\limits_{m = 0}^{2k} {{{\left( { - 1}
\right)}^m}C(m) {x^{2m + 1}}}.
\end{equation}
where $C(0)=C(1)=0$, and for $m\geq2$ the following holds:
$$
C(m) \!=\! \displaystyle\frac{1}{2m+1} -
\displaystyle\frac{4^{m-1}}{3^{m}} \left(
1-8\sum\limits_{i=2}^{m}{\frac{(2i-2)!}{(i-1)!\,i!\,2^{2i-1}}\!\left(\frac{3}{4}\right)^{\!\!i\,}}
\right).
$$

\end{theorem}

\noindent {\bf Proof.} We will prove that the sequence
$\{C(m)\}_{m\in N_{0}}$ is positive, monotonically decreasing and
tends to zero as $m$ tends to infinity.
We will use Lemmas 3~and~4.

\break

\noindent
{\small $$
\begin{array}{lclcl}  
\displaystyle
C(m)
\!&\!\!=\!\!&\!
\displaystyle  \dfrac{1}{{2m + 1}} - \dfrac{3}{2}\dfrac{{(2m - 1)!!}}{{(2m + 2)!!}}\beta_1(m)
\!&\!\!>\!\!&\!
\displaystyle \dfrac{1}{{2m + 1}} - 4 \dfrac{3}{2}\dfrac{{(2m - 1)!!}}{{(2m + 2)!!}} \\[0.75 em]
& & \!&\!\!=\!\!&\!
\displaystyle
\dfrac{1}{{2m + 1}}\left( {1 - 6\dfrac{{(2m + 1)!!}}{{(2m + 2)!!}}}\right).
\end{array}
$$}

\vspace*{-1.0 mm}

\noindent
It is easy to verify that $\displaystyle  6\frac{{(2m + 1)!!}}{{(2m + 2)!!}} < 1$ for $m\geq11$, therefore $C(m)>0$. \\[-0.2 em] 
Let us note that $\displaystyle 0 < C(m) < \frac{1}{{2m + 1}}$, so we can conclude that 
\mbox{$\displaystyle \lim_{m \rightarrow  + \infty}~C(m) = 0$.} \\[-0.2 em]

\noindent
Let us now prove  that $\{C(m)\}_{m \in N_{0}}$ is a monotonically
decreasing sequence.
{\small 
$$
\!\!
\begin{array}{lcl}
C(m+1)  -  C(m)
\!&\!\!=\!\!&\!
\displaystyle \frac{{ - 2}}{{(2m \!+\! 1)(2m \!+\! 3)}} - \frac{3}{2}\frac{{(2m \!+\! 1)!!}}{{(2m \!+\! 4)!!}}{\beta_1}(m\!+\!1) \!+\! \frac{3}{2}\frac{{(2m \!-\! 1)!!}}{{(2m \!+\! 2)!!}}{\beta_1}(m) \\[1.0em]
\!&\!\!=\!\!&\!
\displaystyle  \frac{{ - 2}}{{(2m \!+\! 1)(2m \!+\! 3)}} \!+\! \frac{3}{2}\frac{{(2m\!-\!1)!!}}{{(2m \!+\! 2)!!}}\!\left( {{\beta_1}\!\left(m \right) - \frac{{2m \!+\! 1}}{{2m \!+\! 4}}{\beta_1}\!\left( {m \!+\! 1} \right)} \right) \\[1.0 em]
\!&\!\!<\!\!&\!
\displaystyle  \frac{{ - 2}}{{(2m \!+\! 1)(2m \!+\! 3)}} \!+\! \frac{3}{2}\frac{{(2m\!-\!1)!!}}{{(2m \!+\! 2)!!}}\!\left( {{\beta_1}\!\left( m \right) - \frac{{2m \!+\! 1}}{{2m \!+\! 4}}\!\cdot\!\frac{{m \!+\! 2}}{{m \!+\! \frac{1}{2}}}\!\left( {1 - \frac{3}{2}\frac{1}{{m \!+\! 2}}} \right)\!{\beta_1}\!\left( m \right)} \right) \\[1.0 em]
\!&\!\!=\!\!&\!
\displaystyle \frac{{- 2}}{{(2m \!+\! 1)(2m \!+\! 3)}} \!+\! {\beta_1}\!\left(m\right)\frac{9}{4}\frac{{(2m\!-\!1)!!}}{{\left( {m \!+\! 2} \right)(2m \!+\! 2)!!}}  \\[1.0 em]
\!&\!\!<\!\!&\!
\displaystyle \frac{{ - 2}}{{(2m \!+\! 1)(2m \!+\! 3)}} \!+\! 4\frac{9}{4}\frac{{(2m\!-\!1)!!}}{{\left( {m \!+\! 2} \right)(2m \!+\! 2)!!}} \\ [1.25em]
\!&\!\!=\!\!&\!
\displaystyle \frac{{ - 2}}{{(2m \!+\! 1)(2m \!+\! 3)}}\!\left( {1 - 9\frac{{(2m \!+\! 3)!!}}{{(2m \!+\! 4)!!}}} \right).
\end{array}
$$}

\vspace*{-2.0 mm}

\noindent
It is easy to prove that  $C(m+1) - C(m) < 0$ for $m \ge 8$, i.e. the sequence is monotonically decreasing.
Since $\left\{ {{C }\left( m \right)} \right\}_{m \in N_{0}}\,$is positive for \mbox{$m\!\geq\!2$}, monotonically decreasing
(for \mbox{$m\!\geq\!8$}), and tends to zero, the same holds true for the sequence $\left\{{C(m)x^{2m+1}} \right\}_{m \in N_{0}}$
for~a fixed $\displaystyle x \!\in\! {\big (}0,\sqrt{3}/2{\big ]}$ {\big (}noting that it is decreasing for $m\geq3${\big )}, so we can apply 
{\sc Leibniz}'s theorem for alternating series \cite{Beals_2004}, thus proving the claim of Theorem~\ref{improvement-MS-th-1}:
$$
\sum\limits_{m = 0}^{2k + 1} {{{\left( { - 1} \right)}^m}C(m) {x^{2m + 1}}}
< f(x) <
\sum\limits_{m = 0}^{2k} {{{\left( { - 1} \right)}^m}C(m) {x^{2m + 1}}},  \, k \in N.
$$
\hfill $\Box$

\noindent {\bf Examples.}\\
For $k\!=\!1$ and $x \!\in\! {\big (}0,\sqrt{3}/2{\big ]}$ we get Statement~\ref{MS_statement_1}.\\
For $k\!=\!2$ and $x \!\in\! {\big (}0,\sqrt{3}/2{\big ]}$:
\[\displaystyle
\frac{{{x^5}}}{{180}} - \frac{{13{x^7}}}{{1512}} + \frac{{53{x^9}}}{{5184}} - \frac{{3791{x^{11}}}}{{342144}} < \mbox{\rm arctan}\, x - \frac{{3x}}{{1 + 2\sqrt {1 + {x^2}} }} < \frac{{{x^5}}}{{180}} - \frac{{13{x^7}}}{{1512}} + \frac{{53{x^9}}}{{5184}}.
\]
For $k\!=\!3$ and $x \!\in\! {\big (}0,\sqrt{3}/2{\big ]}$:
\[\displaystyle
\!\!\!
\begin{array}{c}
\displaystyle \frac{{{x^5}}}{{180}}\! -\! \frac{{13{x^7}}}{{1512}}\! +\! \frac{{53{x^9}}}{{5184}}\! -\! \frac{{3791{x^{11}}}}{{342144}}\! +\! \frac{{55801{x^{13}}}}{{4852224}}\! -\! \frac{{130591{x^{15}}}}{{11197440}}\! < \! \mbox{\rm arctan}\, x - \frac{{3x}}{{1 + 2\sqrt {1 + {x^2}} }} \\[1.25 em]
\displaystyle <\frac{{{x^5}}}{{180}}\! -\! \frac{{13{x^7}}}{{1512}}\! +\! \frac{{53{x^9}}}{{5184}}\! -\! \frac{{3791{x^{11}}}}{{342144}}\! +\! \frac{{55801{x^{13}}}}{{4852224}}.
\end{array}\]

\noindent
etc.

\break

\subsection{   Refinements of the inequalities in Statement 2}  

 We propose the following improvement and generalization of Statement \ref{MS_statement_2}:

\begin{theorem}
\label{improvement-MS-th-2} For every $x \in(0,1]$  and    $k \in
N$, it is asserted that
\begin{equation}
\label{improvement2}\displaystyle \frac{3x  +  \mbox{\small
$\displaystyle\sum_{m=2}^{2k+1}$}(-1)^{m} E(m)x^{2m+1}}{1 +
2\sqrt{1+x^2}} <  \mbox{\rm arctan}\,x  <
 \frac{3x + \mbox{\small
$\displaystyle\sum_{m=2}^{2k}$}(-1)^{m} E(m) x^{2m+1}}{1 +
2\sqrt{1+x^2}},
\end{equation}
where
$$
E(m)= \displaystyle
\frac{3}{2m+1}-\sum_{i=0}^{m-1}\frac{(2m-2i-2)!}{2^{(2m-2i
-2)}\,(2i+1)\,(m-i-1)!\,(m-i)!\,}.
$$
\end{theorem}

\noindent
{\bf Examples }\\
For $ x\in (0,  1]$  and   $k  =  1$  we get the inequality
(\ref{MS_2}) from Statement~\ref{MS_statement_2}. \\
For $x \in (0,1]$   and   $k\geq 2$ the inequality
(\ref{improvement2}) refines the inequality (\ref{MS_2}) from
Statement 2 and we have the following new results:
\begin{itemize}
\item
Taking   $k=2$   in  (\ref{improvement2}) gives
{\small
$$
\!\!\!\!\!\!\!\!\!\!\!\!\!\!\!\!\!\!\!\!
\displaystyle \frac{3x \, + \,  \mbox{\small $\displaystyle \frac{1}{60}$}x^5 - \mbox{\small
$\displaystyle \frac{17}{840}$}x^7  +  \mbox{\small $\displaystyle \frac{139}{6720}$}x^9 - \mbox{\small
$\displaystyle \frac{8947}{443520}$}x^{11}}{1+2\sqrt{1+x^{2}}}
<
 \displaystyle\mbox{\rm arctan}\,x
<
\displaystyle \frac{3x \, + \, \mbox{\small $\displaystyle
\frac{1}{60}$}x^5 - \mbox{\small
$\displaystyle \frac{17}{840}$}x^7  + \mbox{\small $\displaystyle \frac{139}{6720}$}x^9}{1+2\sqrt{1+x^{2}}}.
$$
}
\item Taking   $k=3$   in  (\ref{improvement2}) gives

$\begin{array}{l}
\displaystyle \frac{3x \, + \, \mbox{\small $\displaystyle \frac{1}{60}$}x^5 - \mbox{\small
$\displaystyle  \frac{17}{840}$}x^7  +  \mbox{\small $\displaystyle \frac{139}{6720}$}x^9 - \mbox{\small
$\displaystyle \frac{8947}{443520}$}x^{11}
+ \mbox{\small $\displaystyle \frac{89279}{4612608}$}x^{13}
- \mbox{\small $\displaystyle \frac {851677}{46126080}$}x^{15}
}
{1+2\sqrt{1+x^{2}}} \\  [1.5em]
 \, < \,
 \displaystyle\mbox{\rm arctan}\,x
< \,
\displaystyle \frac{3x \, + \, \mbox{\small $\displaystyle
\frac{1}{60}$}x^5 - \mbox{\small
$\displaystyle \frac{17}{840}$}x^7  + \mbox{\small $\displaystyle \frac{139}{6720}$}x^9  - \mbox{\small
$\displaystyle \frac{8947}{443520}$}x^{11}
+ \mbox{\small $\displaystyle \frac{89279}{4612608}$}x^{13}
}
{1+2 \sqrt{1+x^{2}}}.
\end{array}$
\item[] etc.
\end{itemize}


\noindent
{\bf Proof  of Theorem~\ref{improvement-MS-th-2}.}
Based on {\sc Cauchy}'s product of power series (\ref{label_K}) and (\ref{label_A}), the real analytical function
\begin{equation}
f(x) = \left(1 + 2\sqrt{1+x^2}\,\right) \!\cdot\! \mbox{\rm
arctan}\,x - 3x,
\end{equation}
for $x \in  (0,  1]$  has the following power series:
\begin{equation}
f(x) = \sum\limits_{m=2}^{\infty}{(-1)^{m} E(m) \, x^{2m+1}}, ~
\end{equation}
where
\begin{equation}
E(m)= \displaystyle
\frac{3}{2m+1}-\sum_{i=0}^{m-1}\frac{(2m-2i-2)!}{(2i+1)\,(m-i-1)!\,(m-i)!\,2^{(2m-2i
-2)}}.
\end{equation}

We aim to show that sequence $\left\{E(m)\right\}_{m\in N, \;
m\geq2}$  decreases monotonically and that $\displaystyle
\lim_{m\rightarrow +\infty} E(m) =0$. It is easy to verify that
sequence $\left\{E(m)\right\}_{m\in N, \; m\geq 2}$ satisfies the
following recurrence relation:
\begin{equation}
\label{recurence-2}
-2 m E(m)\,+\, (2m+3) E(m+1) \, = \, \displaystyle \frac{1}{2m+1} \,
- \, \frac{(m+1)\, (2m)!}{\,((m+1)!)^{2} \, 4^{m}\,}.
\end{equation}
Consider the sequence $\left\{e(m) \right\}_{m \in N, \; m\geq 2}$ where
$$
\mathit{e}(m) = g(m) \cdot \mathrm{S}(m)
$$
and
$$
S(m) = \displaystyle  {\displaystyle \sum
_{j=1}^{m - 1}} \,{\displaystyle \frac {(2\,j + 2)\mathrm{!}\,
\left( \,((2\, j)\mathrm{!!})^{2}\,(j + 1) - (2\,j + 1 )\mathrm{!}\,\,
\right)\, }{\,2^{(2\,j + 3)}\,(2\,j + 1)\,((j + 1)\mathrm{!})^{2}\,
\,((2\,j)\mathrm{!!})^{2}}}
$$
and
$$
g(m) =  {\displaystyle \frac {m\mathrm{!}\,(m - 1)\mathrm{!} \,2^{(2\,m + 1)}\,}{(2\,m + 1)
\mathrm{!}}}.
$$

It is easy to verify that sequence $\left\{e(m)\right\}_{m\in N, \;
m\geq 2}$ satisfies the   recurrence relation~(\ref{recurence-2}).
Given that  sequences $\left\{E(m) \right\}_{m \in N, \; m\geq 2}$
and $\left\{e(m) \right\}_{m \in N, \; m\geq2}$ agree for $m=2$ and
$m=3,$ we conclude that
\begin{equation}\label{jednakost}
E(m) =e(m),~~ \mbox{for} ~  m\in N, \, m\geq 2.
\end{equation}

We prove that  sequence  $\left\{e(m) \right\}_{m \in N, \; m\geq2}$
is a monotonically decreasing sequence and  $\displaystyle
\lim_{m\rightarrow +\infty} e(m) =0$.

By the principle of mathematical induction, it follows that
$$
(2\,j + 1)\mathrm{!}\,  <  \,(\, (2 \, j)\mathrm{!!}\,)^{2}\, (j +
1)
$$
is true for all $j \in N.$  Therefore  $S(m)>0$  for $m\geq 2$, i.e.
\begin{equation}
\label{pozitivnost}
e(m) > 0,  ~~ \mbox{for } ~  m\geq 2.
\end{equation}
 To prove that   $\left\{e_{m}\right\}_{m \in N,\, m\geq 2}$   is a
monotonically decreasing   sequence, let us  use the following
notation:
$$
S(m) = \displaystyle \sum_{j=1}^{m - 1} h(j)
$$
where
$$
\displaystyle h(j) = {\displaystyle \frac {(2\,j + 2)\mathrm{!}\,
\left( ((2\, j)\mathrm{!!})^{2}\,(j + 1) - (2\, j + 1 )\mathrm{!}\,
\right) }{2^{(2\,j + 3)} \, (2\,j + 1)\,((j + 1)\mathrm{!})^{2}\,
\,((2\,j)\mathrm{!!})^{2}}}.
$$

\medskip
Consider the  following equivalences for  $m\geq 2:$
\begin{equation}
\begin{array}{rcl}
\displaystyle  \frac{e(m+1)}{e(m)}< 1
\!&\!\Longleftrightarrow\!&\!
\displaystyle   \frac{g(m+1)\, S(m+1)}{g(m) \, S(m)} < 1 \\  [1.2em]
\!&\!\Longleftrightarrow\!&\!
\displaystyle  \frac{g(m+1)}{g(m)} \cdot  \frac{S(m)+h(m)}{S(m)}  < 1  \\  [1.2em]
\!&\!\Longleftrightarrow\!&\!
\displaystyle \frac{2m}{2m+3} \cdot \frac{S(m)+h(m)}{S(m)}  < 1  \\ [1.2em]
\!&\!\Longleftrightarrow\!&\!
\displaystyle   2\, m \, h(m) < 3\, S(m).
\end{array}
\end{equation}

\medskip
Consider the last inequality.  It is easy to verify that it is true for $m=2$. Observing that
$$
\begin{array}{rcl}\displaystyle  3 \, S(m+1)
\!&\!=\!&\! 3 \left( \, S(m) \,  + \, h(m) \, \right) \\[1.0 ex]
\!&\!=\!&\! 3 \, S(m) \, + \, 3 \,h(m)
\end{array}$$
and using the induction hypothesis {\large $(\!$} $3 \, S(m) >  2 \, m \, h(m)$  {\large $\!)$}
for some positive integer $m \geq 2$, we conclude that
$$
3 \, S(m+1) > (2m+3) h(m) > 2 (m+1) h(m).
$$
Therefore, by the principle of mathematical induction, the inequality
\begin{equation}
\displaystyle   2\, m \, h(m)  < 3\, S(m)
\end{equation}
is true for  $m\geq 2$  i.e.
\begin{equation}\label{monotonost}
 \displaystyle  \frac{e(m+1)}{e(m)}< 1, \quad \mbox{for $m \geq  2$.}
\end{equation}

Let us  further consider the   positive addend  of $S(m)$,  i.e.
$$
\mathit{S}^{+}(m) = {\displaystyle \sum _{j=1}^{m - 1}}
\,{\displaystyle \frac {(2\,j + 2)\mathrm{!}\,\,((2 \, j)\mathrm{!!}
\,)^{2}\,(j + 1)}{\,  2^{(2\,j + 3)}  \, (2\,j + 1)\,((j +
1)\mathrm{!})^{2}\, (\,(2 \,j)\mathrm{!!}\,)^{2}}} =
 {\displaystyle \sum _{j=1}^{m - 1}}
\,{\displaystyle \frac {(2\,j)\mathrm{!}\,}{4 \,((2 \,j
)\mathrm{!!})^{2}}}.
$$

\noindent By the principle of mathematical induction,  it  follows
that

\[ \mathit{S}^{+}(m) =
{\displaystyle \frac{m  \,(2\,m)\mathrm{!}\, }{\,2 \,\,
((2\,m)\mathrm{!!})^{2}\,} \, }   -  {\displaystyle  \frac{1}{4}  }.
\]

Finally, given that for $m \geq 2$
$$ 0 \,< \, e(m) \,  \leq  \,  g(m) \cdot \mathit{S}^{+}(m) \,  = \,   \displaystyle
\frac{1}{2m+1}$$ we have $\displaystyle \lim_{m \rightarrow +\infty}
e(m) = 0.$

~\\
Finally, based on (\ref{jednakost}) we  conclude
that $\left\{ {E\left( m \right)} \right\}_{m\in N, m\geq 2}$ is a positive
 monotonically decreasing sequence, and that it tends to
zero. The same holds true for the sequence $\left\{{E(m)x^{2m+1}}
\right\}_{m\in N, m\geq 2}$ for a fixed $x \in (0,1]$  so we can
apply {\sc Leibniz}'s theorem for alternating series \cite{Beals_2004},
thus proving the claim of Theorem~\ref{improvement-MS-th-2}.

\hfill $\Box$

\bigskip
\subsection{   Refinements of the inequalities in Statement 3}

We propose the following improvement and generalization of Statement
3:

\begin{theorem}
\label{improvement-MS-th-4} For every $x \in(0,1]$  and    $k \in
N$, it is asserted that
\begin{equation} \label{improvement4}
\displaystyle \sum_{m=1}^{2k-1}(-1)^{m} C(m)x^{2m+1}
\!<\!
\mbox{\rm arctan}\,x - \frac{2x }{\,1 + \sqrt{1+x^{2}}\,}
\!<\!
\displaystyle\sum_{m=1}^{2k}(-1)^{m} C(m) x^{2m+1},
\end{equation}
where
$$
C(m)= \displaystyle \frac{1}{2m+1}- \frac{(\,2 m-1\,)!!}{(m+1)! \,
2^{m} }.
$$
\end{theorem}

\noindent  {\bf Examples }\\
For $ x\in (0,  1]$  and   $k  =  1$  we get the inequality
(\ref{MS_4})  from Statement~\ref{MS_statement_4}. \\
For $x \in (0,1]$   and   $k\geq 2$ the inequality
(\ref{improvement4}) from Theorem \ref{improvement-MS-th-4}  refines
the inequality (\ref{MS_4}) from Statement~\ref{MS_statement_4} and
we have the following new results:
\begin{itemize}
\item
Taking   $k=2$   in  (\ref{improvement4}) gives
$$
\!\!\!\!\!\!\!\!\!\!\!\!\!\!\!\!
\mbox{\small $\displaystyle -\frac{1}{12}$}x^3 + \mbox{\small
$\displaystyle \frac{3}{40}$}x^5 - \mbox{\small $\displaystyle
\frac{29}{448}$} x^{7}
<
 \displaystyle\mbox{\rm arctan}\,x -
\frac{2x}{1 + \sqrt{1-x^2}}
<
\mbox{\small $\displaystyle
-\frac{1}{12}$}x^3 + \mbox{\small $\displaystyle \frac{3}{40}$}x^5 -
\mbox{\small $\displaystyle \frac{29}{448}$} x^{7}+ \mbox{\small
$\displaystyle \frac{65}{1152}$} x^{9}.
$$
\item Taking   $k=3$   in  (\ref{improvement4}) gives
$$
\!\!\!\!
\begin{array}{ll}
\mbox{\small $\displaystyle -\frac{1}{12}$}x^3 + \mbox{\small
$\displaystyle \frac{3}{40}$}x^5 - \mbox{\small $\displaystyle
\frac{29}{448}$} x^{7} & \!\!\!\!\!\! + \, \mbox{\small $\displaystyle
\frac{65}{1152}$} x^{9} -  \mbox{\small $\displaystyle
\frac{281}{5632}$} x^{11} <   \displaystyle\mbox{\rm arctan}\,x -
\frac{2x}{\,1 + \sqrt{1-x^2}\, }  \\  [1.5em]
&
\!\!\!\!\!\!\!\!\!\!\!\!\!\!\!\!\!\!\!\!\!\!\!\!
\!\!\!\!\!\!\!\!\!\!\!\!\!\!\!\!\!\!\!\!\!\!\!\!
<\mbox{\small $\displaystyle
-\frac{1}{12}$}x^3 + \mbox{\small $\displaystyle \frac{3}{40}$}x^5 -
\mbox{\small $\displaystyle \frac{29}{448}$} x^{7}+ \mbox{\small
$\displaystyle \frac{65}{1152}$} x^{9}  -  \mbox{\small
$\displaystyle \frac{281}{5632}$} x^{11}  +  \mbox{\small
$\displaystyle \frac{595}{13312}$} x^{13}.
\end{array}
$$
\item[] etc.
\end{itemize}

\noindent
{\bf Proof  of Theorem~\ref{improvement-MS-th-4}.}

For $x \in (0,1]$ the following power series expansion holds:
$$\begin{array}{rcl}
\displaystyle \mbox{\rm arctan}\,x \, -\, \frac{2x}{1+\sqrt{1+x^2}}
\!&\!\!=\!\!&\!
\displaystyle \mbox{\rm arctan}\,x \, +\,  {\displaystyle \frac {2\,(1 - \sqrt{x^{2} + 1})}{x}} \\ [1.5em]
\!&\!\!=\!\!&\!
\displaystyle\sum_{m=1}^{\infty}{(-1)^{m} C(m) \, x^{2m+1}}
\end{array}
$$
where
\begin{equation}
\label{C(m)_T4} C(m) = {\displaystyle \frac {1}{2\,m + 1}}
 - {\displaystyle \frac {(2\,m - 1)\mathrm{!!}}{(m
 + 1)\mathrm{!}\,2^{m}}}.
\end{equation}

We  prove that the sequence $\left\{C(m)\right\}_{m\in N}$ is
positive, monotonically decreasing and  tends to zero as $m$ tends
to infinity.

It is easy  to verify that $ \displaystyle (2m + 2)!! \,  >  \, \displaystyle 2 \,  (2m + 1)!! $
is true for $m \in N$.  Thus,  the following equivalences  hold true for every $\,m\in N$
$$
\begin{array}{lcl}
C(m)   > 0
\!&\!\!\Longleftrightarrow\!&\!\!
\displaystyle \frac{1}{2m+1} > \displaystyle  \frac{(\,2 m-1\,)!!}{(m+1)! \, 2^{m} }        \\[2.0 ex]
 & \Longleftrightarrow & \displaystyle (2m + 2)!! \,  >  \, \displaystyle 2 \,  (2m + 1)!!,
\end{array}
$$
and  we conclude that  $C(m)>0$  for every $m \in N$.

Let us  now prove that $\left\{ C(m) \right\}_{m\in N}$ is a
monotonically decreasing sequence.  We have:
\[\begin{array}{lcl}
C(m)  -  C(m+1)   > 0
\!&\!\!\Longleftrightarrow\!\!&\!
\displaystyle \frac{2}{\,(2m+1)\,(2m+3)\,} - 3 \frac{ \,(2m -1)!!}{\,(2m +2)!! \, (m+2)} > 0  \\ [2.0 ex]
\!&\!\!\Longleftrightarrow\!\!&\!
\displaystyle (2m+4)!! \, > \, 3 \, (2m+3)!!.
\end{array}\]

As it  is easy to show (by the principle of mathematical induction)
that  the last inequality holds true  for $m \in N$, we may conclude
that $\left\{ C(m) \right\}_{m\in N}$ is a monotonically decreasing
sequence.

Finally, as  $0 < C(m)  < \mbox{\small $\displaystyle \frac{1}{{2m + 1}}$}$,  we
conclude that $\displaystyle \lim_{m \rightarrow +\infty}  C(m)  =  0$.

Since $\left\{ {C\left( m \right)} \right\}_{m\in N}$ is a positive
monotonically decreasing sequence, and it tends to zero, the same
holds true for the sequence $\left\{{C(m)x^{2m+1}}\right\}_{\!m \in
N}$ for a fixed $x \!\in\! (0,1]$. So we can apply {\sc Leibniz}'s
theorem for alternating series \cite{Beals_2004} and thus prove the
claim of Theorem~\ref{improvement-MS-th-4}.

\hfill $\Box$

\section{Conclusion}

In Theorems 1, 2 and 3 of this paper we proved some new inequalities~related 
to~{\sc Shafer}'s inequality for the arctangent function.$\;$These inequalities 
re\-present sharpening and generalization of the inequalities given  in \cite{MS_2014} (Theorems 1,~2~and~4).

\smallskip
Let us mention that it is possible to prove the inequality (\ref{Mara22}), for any fixed $k \!\in\! N$ and 
$x \!\in\! {\big (}0, \frac{\sqrt{3}}{2}{\big ]}$, by substituting $x \!=\! \tan t$ 
for $t \!\in\! {\big (}0, \mbox{\rm arctan} \frac{\sqrt{3}}{2} {\big ]}$ using the algorithms and
methods${}^{\ast\ast)}$ developed in \cite{MM_2016}~and~\cite{LMM_2017}. 
Also, the inequalities (\ref{improvement2}) and (\ref{improvement4}) for any fixed 
$k \in N$ and $x \!\in\! (0, 1]$  can be  proved  by substituting \mbox{$x \!=\! \tan t$} 
for $t \!\in\! {\big (}0, \frac{\pi}{4}{\big ]}$ using the algorithms 
and methods${}^{\ast\ast)}$\footnote{\!\!\!\!\!\!\!\!\!\!${}^{\ast\ast)}\,$See 
also \cite{Malesevic_Lutovac_Banjac_2016} and \cite{Malesevic_Banjac_Jovovic_2017}.} 
developed in \cite{MM_2016}~and~\cite{LMM_2017}. 

\medskip
\noindent
\textbf{Acknowledgement.} Research of the first, second
and third author was supported in part by the Serbian Ministry of
Education, Science and Technological Development, under  Projects ON
174032, III 44006 and ON 174033 and TR 32023, respectively.

\medskip
\noindent
\textbf{Competing Interests.} The authors would like to
state that they do not have any competing interests in the subject
of this research.

\medskip
\noindent
\textbf{Author's Contributions.} All the authors participated in every phase of the research conducted for this paper.

\break

\end{document}